\documentclass{article}
\usepackage{arxiv}
\usepackage{amsmath,amssymb,amsfonts} \usepackage{orcidlink}
\usepackage{wrapfig}
\usepackage{subcaption}  
\usepackage{amsfonts}
\usepackage{amsopn}
\usepackage{amsmath,amssymb}
\usepackage{graphicx,upref,mmap}
\usepackage{epstopdf}
\usepackage{algorithmic}
\usepackage{amsmath,amssymb,amsfonts} \usepackage{orcidlink}
\usepackage{wrapfig}
\usepackage{xcolor,color}
\usepackage{subcaption}  
\usepackage{tikz-cd}
\usetikzlibrary{positioning,arrows.meta}

\graphicspath{{./},{./figures/}}

\newcommand{\fn}{{\mathrm{fn}}}

\newtheorem{thm}{Theorem}
\newtheorem{lemma}{Lemma}
\newtheorem{prop}{Proposition}
\newtheorem{remark}{Remark}
\def\spacingset#1{\def\baselinestretch{#1}\small\normalsize}
\spacingset{1}

\renewcommand{\fn}{{\mathrm{fn}}}

\newcommand{\factor}{\Phi}
\DeclareMathOperator{\diag}{diag}
\DeclareMathOperator{\rank}{rank}
\DeclareMathOperator{\trace}{trace}
\newcommand{\Sim}{\text{Sym}^+}

\begin{document}

\title{
The factorization of matrices into products of 
positive definite factors
}

\author{Mahmoud Abdelgalil\orcidlink{https://orcid.org/0000-0003-1932-5115}\thanks{Electrical and Computer Engineering, University of California, San Diego, La Jolla, CA, USA, mabdelgalil@ucsd.edu.}
\hspace*{3pt}and Tryphon T. Georgiou\orcidlink{https://orcid.org/0000-0003-0012-5447} 
\thanks{Mechanical and Aerospace Engineering, University of California, Irvine, Irvine, CA, USA, tryphon@uci.edu.\newline 
\hspace*{10pt}{Research supported by NSF under ECCS-2347357, AFOSR under FA9550-24-1-0278, and ARO under W911NF-22-1-0292.}
}}

\maketitle

\begin{abstract}
Positive-definite matrices materialize as state transition matrices of linear time-invariant gradient flows, and the composition of such materializes as the state transition after successive steps where the driving potential is suitably adjusted. Thus, factoring an arbitrary matrix (with positive determinant) into a product of positive-definite ones provides the needed schedule for a time-varying potential to have a desired effect. The present work provides a detailed analysis of this factorization problem by lifting it into a sequence of Monge-Kantorovich transportation steps on Gaussian distributions and studying the induced holonomy of the optimal transportation problem. From this vantage point we determine the minimal number of positive-definite factors that have a desired effect on the spectrum of the product, e.g., ensure specified eigenvalues or being a rotation matrix. Our approach is computational and allows to identify the needed number of factors as well as trade off their conditioning number with their actual number. 
 
\end{abstract}

%

\section{Introduction}
Whereas it is well known that the product $\Phi_1\Phi_2$ of two symmetric positive-definite\footnote{Henceforth, `symmetric positive-definite' will be abbreviate to  `positive' and denoted by $>0$.} matrices $\Phi_1$ and $\Phi_2$ has positive eigenvalues, $\Phi_1\Phi_2$ being similar to the positive-definite matrix $\Phi_1^{1/2}\Phi_2\Phi_1^{1/2}$, the situation is more nuanced for products of three symmetric matrices or more. Indeed, while positivity of the eigenvalues fully characterizes diagonalizable matrices as products of two positive ones, and the factors are straightforward to obtain\footnote{If $M=TQT^{-1}$, with $T$ invertible and a diagonal $Q>0$, then $M=\Phi_1\Phi_2$ for $\Phi_1=TT^*>0$ and $\Phi_2=T^{-*}QT^{-1}>0$. Throughout, $*$ denotes complex conjugation and transposition.},
no \emph{simple} construction exists in the literature for matrices with negative or complex eigenvalues, in which case the least number of factors exceeds two \cite[Theorem 3]{ballantine1970products}.
It is however known that, regardless of the dimension, a product of five positive factors is sufficient to factor an arbitrary real square matrix with positive determinant, as shown by Ballantine in a series of publications \cite{ballantine1967products,ballantine1968aproducts,ballantine1968bproducts,ballantine1970products}. Our goal is to elucidate such `Ballantine factorizations' via a principled constructive method for obtaining the required positive factors. 
Our approach is geometric, viewing the factorization problem as a succession of gradient flows that translocate constituent particles of Gaussian distributions. By lifting the problem to one that concerns the holonomy of a transportation problem \cite{abdelgalil2025holonomy}, it provides a bird's eye view on the freedom in the selection of factors, a systematic method to select such, and the ability to control the conditioning number of the factors.

Specifically, we build on the authors' recent work that introduces a sub-Riemannian structure in Monge-Kantorovich optimal transport \cite{abdelgalil2025holonomy}, and quantify the holonomy accrued by 
a positive map
\[
\Phi\;:\; x\mapsto \Phi x
\]
with $x\in\mathbb R^n$ and $\Phi=\Phi^\top \succ0$ on Gaussian distributions. This allows to conveniently book-keep the cumulative effect of a product of positive matrices that effect steps in successive gradient flows. The result of the approach is to provide a systematic as well as `visual' computational procedure that can be used to answer questions that pertain to how many factors are needed and to how much freedom one has in choosing those factors (that are highly non-unique).

The control theoretic significance of factorizations into positive factors was brought to light by the authors in a recent contribution \cite{abdelgalil2024collective} where it is proven that the bi-linear dynamical system
\begin{equation}\label{eq:basic}
    \cfrac{d}{dt}\Phi(t) = (A+BK(t))\Phi(t),
\end{equation}
with $\Phi(t)\in\mathbb R^{n\times n}$, regulated via a time-varying state-feedback gain matrix $K(t)\in {\mathbb R}^{m\times n}$, is strongly controllable if and only if the pair $(A,B)\in{\mathbb R}^{n\times n}\times {\mathbb R}^{n\times m}$ satisfies the Kalman rank condition
\begin{equation}
    \rank\left(\begin{bmatrix}
        B,& AB,& \ldots & A^{n-1}B
    \end{bmatrix}\right)=n,
\end{equation}
with $n,m$ positive integers.
Conditions for the controllability of
\eqref{eq:basic}
were claimed by Roger Brockett in his influential paper on the Liouville equation \cite{brockett2007optimal}, though the approach in \cite{brockett2007optimal} fell short of establishing strong controllability (see \cite{abdelgalil2024collective} for a detailed discussion).
The crux of the matter, for establishing strong controllability in \cite{abdelgalil2024collective},
was to produce a finite number of steps as control primitives for steering \eqref{eq:basic} through
intermediary values of $\Phi(t_i)$ from the identity $\Phi(t_0)=I$ at time $t=t_0$ to the specified value $\Phi(t_\fn)$ at a terminal time $t_\fn$. Moreover, these steps must be such that they can be traversed arbitrarily fast via a suitable choice of a control law $K(\cdot)$; this was accomplished by ensuring gradient flow in each step, which in turn amounts to factoring $\Phi(t_\fn)$ as a product of symmetric matrices with a known bound on the number of such factors.
As far as we are aware, \cite{abdelgalil2024collective} provides the first articulation of Ballantine-type  factorizations as a control-theoretic device. It shows that gradient-driven flows can be traversed in arbitrarily small time without inducing shocks, which in turn yields strong controllability when applicable.

In the sequel, in Section \ref{sec:angular} we highlight the conceptual framework of the approach along with reviewing some basic concepts of the theory of optimal transport, specialized to Gaussian distributions. In particular, we highlight a concept of holonomy incurred by optimal transport cycles. Section \ref{sec:3} builds on this concept of holonomy to quantify rotation accrued along succession of optimal transport maps. In Section \ref{sec:4} we present the main results, first characterizing invertible $2\times 2$ matrices that admit factorization into a product of $2$, $3$, $4$, and $5$ (five being the maximal number of factors needed), positive factors. Section \ref{sec:5} we extend the results to $n\times n$ invertible matrices, and conclude with remarks in Secction \ref{sec:remarks}.






\section{Holonomy in the transport of Gaussian distributions}\label{sec:angular}

We consider centered Gaussian probability distributions in $\mathbb R^n$ with covariance $\Sigma_i$, namely,
\[
\mathcal G_i(x)=(2\pi)^{-n/2}\det(\Sigma_i)^{-1/2}\exp(-x^\top\Sigma_i^{-1}x/2), \mbox{ for }x\in\mathbb R^n,
\]
for $i\in\{1,2,\ldots,k\}$; at times, we also use the notation $N(0,
\Sigma)$ for a centered Gaussian with covariance $\Sigma$. The covariance matrices $\Sigma_i$ are assumed to be positive definite throughout. Scaling is immaterial for our purposes, and thereby we assume that $\det(\Sigma_i)=1$. We are interested in the change of the relative position of constituent particles with respect to the principal axes of distributions, effected by optimally transferring one distribution into another. 
The transportation is optimal in the Monge-Kantorovich sense, i.e., it minimizes the Wasserstein distance between the two distributions and is effected by a positive map. Change in the relative position of constitutive particles amounts to a holonomy accrued by the transportation map in the sense of \cite{abdelgalil2025holonomy}. This we explain next.

\subsection{Optimal Transport}
Consider a proto-typical particle at location $x_1\in\mathbb R^n$, distributed according to $\mathcal G_1$ and transported into
\[
x_2 = \Phi_{21}x_1,
\]
with $\Phi_{21}\in\mathbb R^{n\times n}$, so that $x_2$ is now distributed according to $\mathcal G_2$ with covariance $\Sigma_2=\Phi_{21}\Sigma_1\Phi_{21}^\top$. Given the two covariance matrices $\Sigma_1$ and $\Sigma_2$, the choice of $\Phi_{21}$ is, evidently, non-unique. However, the choice that minimizes the quadratic cost
\[
\mathbb E_{\mathcal G_1}\| x_2-x_1\|^2,
\]
subject to $x_2$ being distributed according to $\mathcal G_2$, is unique; the optimal choice for $\Phi_{21}$ is symmetric and positive definite, and given by either of the following two equivalent expressions (see \cite{dowson1982frechet,olkin1982distance})
\begin{align}\label{eq:knott}
\Phi_{21}^\star&=\Sigma_1^{-1/2}\left(
\Sigma_1^{1/2}
\Sigma_2
\Sigma_1^{1/2}
\right)^{1/2}
\Sigma_1^{-1/2}\\\nonumber
&=\Sigma_2^{1/2}\left(
\Sigma_2^{1/2}
\Sigma_1
\Sigma_2^{1/2}
\right)^{-1/2}
\Sigma_2^{1/2}.
\end{align}
Also note that  $\Phi_{12}^\star\Phi_{21}^\star=I$. The solution of this optimal transportation problem, taking $x_1$ to $x_2$ while minimizing the aforementioned quadratic cost, is completely symmetric with respect to transporting $x_1$ to $x_2$, or the other way around. That is, $\Phi_{12}^\star$ is the optimal transport map that minimizes $\mathbb E_{\mathcal G_2}\| x_1-x_2\|^2$ when $x_1$ is sought as a function of $x_2$, and required to be distributed according to $\mathcal G_1$. 

The same solution as above can also be obtained by seeking an optimal joint distribution $\mathcal G$ of $x_1$ and $x_2$ that minimizes
    $\mathbb E_{\mathcal G} \|x_2-x_1\|^2$,
while having the specified marginals $\mathcal G_1$ and $\mathcal G_2$.
This so-called `coupling' is a centered Gaussian
 with covariance
\begin{equation}\label{eq:coupling}
    \Sigma= \begin{bmatrix}
       \Sigma_1 & \Sigma_{12}\\
        \Sigma_{12}^\top & \Sigma_2
    \end{bmatrix}, \mbox{ with }\Sigma_{12}:=\mathbb E\{x_1x_2^\top\}=\Sigma_1\Phi_{21}^\star=\Phi_{12}^\star\Sigma_2.
\end{equation}

Either way, the minimal quadratic cost constitutes a metric in the space of Gaussian distribution, the Wasserstein-$2$ metric of the Monge-Kantorovich theory \cite{villani2003topics,villani2009optimal}, and this is
\[
W_2(N(0,\Sigma_1),N(0,\Sigma_2))^2=\trace\left(\Sigma_1+\Sigma_2-2(\Sigma_1^{1/2}\Sigma_2\Sigma_1^{1/2}\right).
\]

\subsection{Holonomy of transport}
We view $x_1\in\mathbb R^n$ as the position of a particle, distributed according to $\mathcal G_1$, and are interested in the relative angular displacement of all such $\mathcal G_1$-particles as effected by the optimal transportation map $\Phi_{21}^\star$.
The central idea of our approach stems from realizing that optimal transport amounts to gradient flow, and a composition of such presents a natural framework to study factorization of maps as sought herein.
The holonomy of optimal transport, quantifying distance from the identity when completing closed cycles in the Wasserstein space, was recently introduced by the authors in \cite{abdelgalil2024sub}. We will associate a holonomy with any leg of a path by appealing to certain canonical triangular cycles in the Wasserstein space; in other words, by registering the translocation of the particles relative to a canonical frame. We outline this next.

We identify a non-degenerate centered Gaussian distributions in $\mathbb R^n$ with the corresponding covariance matrix $\Sigma \in \Sim(n)$, where $\Sim(n)$ is the space of $n\times n$ real symmetric positive definite matrices. The smooth submersion:
\begin{align}
    \pi\;:\;\text{GL}^+(n)\to \text{Sym}^+(n)\;:\;
    \factor\mapsto \factor \factor^\top
\end{align}
where $\text{GL}^+(n)$ is the identity component  of the general linear group, defines a smooth fiber-bundle with total space $\text{GL}^+(n)$ over the base space $\Sim(n)$. For any $\Sigma\in\text{Sym}^+(n)$, the fiber at $\Sigma$ is the submanifold:
\begin{align}
    \pi^{-1}(\Sigma) = \{\factor\in\text{GL}^+(n)~|~\pi(\factor) = \Sigma\}.
\end{align}
In particular, the fiber $\pi^{-1}(I)$ coincides with the subgroup
\begin{align}
    \text{SO}(n)=\{\Theta\in\text{GL}^+(n)~|~ \Theta\Theta^\top = I\},
\end{align}
which is the \emph{isotropy group} of the covariance $I\in\Sim(n)$ under the left action of $\text{GL}^+(n)$ on $\Sim(n)$ defined by: $\factor\cdot\Sigma:=\factor\Sigma \factor^\top$. As a (closed) subgroup of $\text{GL}^+(n)$, $\text{SO}(n)$ has a free right action on $\text{GL}^+(n)$ via matrix multiplication:
\begin{align}
    R_\Theta\;:\;\text{GL}^+(n)\to \text{GL}^+(n)\;:\; \factor \mapsto \factor\Theta.
\end{align}
In addition, the right action of $\text{SO}(n)$ on $\text{GL}^+(n)$ preserves the fibers of the bundle $\pi$ since, if $\factor\in\pi^{-1}(\Sigma)$ and $\Theta\in\text{SO}(n)$, then
\begin{align}
    \pi\circ R_\Theta(\factor) = \factor\Theta \Theta^\intercal \factor^\intercal = \factor \factor^\intercal =  \Sigma.
\end{align}
Finally, it can be shown that the action of $\text{SO}(n)$ on each fiber is \emph{transitive} \cite{modin2016geometry}, i.e. if $\factor$ and $\tilde{\factor}$ belong to the same fiber $\pi^{-1}(\Sigma)$, then there exists $\Theta\in \text{SO}(n,\Sigma_{\text{ref}})$ such that $\tilde{\factor}=\factor\Theta$. Since the action of $\text{SO}(n)$ on the fibers is both free and transitive, the following proposition is immediate \cite{montgomery2002tour,kobayashi1996foundations}.
\begin{prop}\normalfont
    $\pi:\text{GL}^+(n)\rightarrow\text{Sym}^+(n)$ is a principal $\text{SO}(n)$-bundle.
\end{prop}

It is not difficult to see that the principal bundle $\pi$ is trivial, i.e., it admits a global section. Indeed, the map 
\begin{align}
    \sigma:\text{Sym}^+(n)\rightarrow\text{GL}^+(n)~~~:~~~\Sigma\mapsto \Sigma^{\frac{1}{2}}, 
\end{align}
satisfies $\pi\circ\sigma(\Sigma) = \Sigma$, for all $\Sigma\in \text{Sym}^+(n)$.

The significance of the above stems from the fact that the identity can serve as a point to `register' paths, by linking their starting and ending point to the identity. Thereby, the isotropy group can fully characterize the translocation of $\mathcal G$-particles as they traverse any given segment in the Wasserstein space.

\subsection{Angular displacement of segments}

Thus, the concept of holonomy allows us quantify the angular displacement that a linear transportation map (symmetric matrix) imparts on particles of a Gaussian distribution along paths, i.e., it provides a principled way to `register' particles with respect to a common reference. To do this, we consider the holonomy accrued over a three-leg cycle,
\[
v_1\mapsto x_1=\Sigma_1^{1/2}v_1\mapsto x_2=\Phi_{21}x_1\mapsto v_2=\Sigma_2^{-1/2}x_2,
\]
that transports the reference distribution $N(0,I)$ back to itself. The sequence 
\begin{equation}\label{eq:sequence}
    N(0,I)\to N(0,\Sigma_1)\to N(0,\Sigma_2)\to N(0,I)
\end{equation}
traces the vertices of a triangle in the space of Gaussian distributions, cf.\ \cite{abdelgalil2025holonomy}. 

The holonomy of the cycle is quantified by a element of the isotropy group $SO(n)$, that constitutes an angular displacement of particles. Our interest is in studying how the accrued holonomy depends on the choices of $\Sigma_1$ and $\Sigma_2$, since we may now associate this with the optimal transportation segment between the two. The angular displacement is quantified by the rotation matrix $\Theta_{21}$ that takes
the normalized vector $v_1:=\Sigma_1^{-1/2}x_1$ into $v_2:=\Sigma_2^{-1/2}x_2$, schematically,
\begin{equation}\label{eq:tikz}
\begin{tikzcd}
x_1 \arrow[r,"\Phi_{21}"] \arrow[d, "\Sigma_1^{-1/2}"'] & x_2 \arrow[d, "\Sigma_2^{-1/2}"] \\
v_1 \arrow[r,"\Theta_{21}"] & v_2
\end{tikzcd}
\end{equation}
where\footnote{We use the same symbol for linear maps and the corresponding matrix representation.}
$\Theta_{21}\;:\; v_1 \mapsto v_2 = \Theta_{21}v_1$,
    giving
    \begin{align}
   \label{eq:R}
\Theta_{21}&=\Sigma_2^{-1/2}\Sigma_1^{-1/2}\left(
\Sigma_1^{1/2}
\Sigma_2
\Sigma_1^{1/2}
\right)^{1/2}=\left(
\Sigma_2^{1/2}
\Sigma_1
\Sigma_2^{1/2}
\right)^{-1/2}
\Sigma_2^{1/2}\Sigma_1^{1/2}.
\end{align}

The amount of rotation, in $\mathbb R^n$, can be quantified by any unitarily invariant norm on $\Theta_{21}$, though this is not needed as we will next focus on rotations in $\mathbb R^2$. The schematic in \eqref{eq:tikz} is a key concept in that it provides a canonical `registration' of particles at locations $x_1$ and $x_2$, to be compared relative to their relative position with respect to (e.g., the principal axes of) the respective Gaussian distributions. 

%

\section{Rotational dispacement in $\mathbb R^2$}\label{sec:3}

We specialize to the case of $2$-dimensional Gaussian distributions so as to derive explicit results on the holonomy that is accrued by suitably selected directions of transport. As we will see, quantifying the holonomy accrued in the $2$-dimensional case, is key to understanding how to factor $2\times 2$ as well as $n\times n$ matrices into positive factors  -- the latter reduces to the $2\times 2$ case.

\subsection{Holonomy of Gaussian triangles in $\mathbb R^2$}
Our task is to assess the holonomy that is accrued traversing the cycle in \eqref{eq:sequence}, for suitable choices of $\Sigma_1$ and $\Sigma_2$.

Since we now work in $\mathbb R^2$, with unit determinants, there are only three parameters that specify the relative positioning and shape of the two distribution $N(0,\Sigma_1)$ and $N(0,\Sigma_2)$. In what follows, we use the symbol $\Theta$ for orthogonal matrices. The indexing may indicate
the corresponding angle, as in $\Theta_\psi$ and $\Theta_\phi$, but may also mark with numerals as in $\Theta_{21}$ a pair of corresponding covariances, being the holonomy to a particular pair, as in \eqref{eq:R}. Taking as our starting coordinate frame the principal axes of $\Sigma_1$, we have that
\begin{align}
\Sigma_1=\begin{bmatrix}\lambda_1&0\\0&1/\lambda_1
\end{bmatrix}, \mbox{ and }
\Sigma_2= \Theta_\psi \begin{bmatrix}\lambda_2&0\\0&1/\lambda_2
\end{bmatrix}\Theta_\psi^\top,\label{eq:Sigma12}
\end{align}
where 
\[
\Theta_\psi:=\begin{bmatrix}
    \cos(\psi)& -\sin(\psi)\\\sin(\psi)& \phantom{-}\cos(\psi)
\end{bmatrix}.
\]
Evidently, this representation is unique for $\psi\in(-\pi/2,\pi/2)$, and, without loss of generality, we may take $\lambda_1\geq 1$.
Substituting in \eqref{eq:R}, the accrued holonomy is
\begin{align}\label{eq:key}
\Theta_{\phi}:=&\Theta_\psi \begin{bmatrix}1/\sqrt{\lambda_2}&0\\0&\sqrt{\lambda_2}
\end{bmatrix}\Theta_\psi^\top\begin{bmatrix}1/\sqrt{\lambda_1}&0\\0&\sqrt{\lambda_1}\end{bmatrix}\times\\
\nonumber&\times
\left(
\begin{bmatrix}\sqrt{\lambda_1}&0\\0&1/\sqrt{\lambda_1}
\end{bmatrix}\Theta_\psi\begin{bmatrix}{\lambda_2}&0\\0&1/{\lambda_2}
\end{bmatrix}\Theta_\psi^\top\begin{bmatrix}\sqrt{\lambda_1}&0\\0&1/\sqrt{\lambda_1}
\end{bmatrix}
\right)^{\frac12},
\end{align}
where, the angle $\phi\in(-\frac{\pi}{2},\frac{\pi}{2})$ is a function of $\lambda_1$, $\lambda_2$ and $\psi$. 

If we now denote the $(i,j)$-entry of $\Theta_\phi$ by $(\Theta_\phi)_{ij}$, we are interested in
\begin{equation}\label{eq:phimax}
\phi(\lambda_1,\lambda_2,\psi):=\tan^{-1}\left( \frac{(\Theta_\phi)_{21}}{(\Theta_\phi)_{11}}\right).
\end{equation}
We make the following claims:
\begin{lemma}\label{lemma:3.1} For any $\lambda,\lambda_1>1$, $\lambda_2>0$, and $\psi\in(-\pi/2,\pi/2)$, the following hold:
\begin{itemize}
    \item[i)]  $\phi(\lambda_1,\lambda_2,\psi)=\phi(\lambda_2,\lambda_1,\psi)$,
   \item[ii)]
$\phi(\lambda,\lambda,\psi)$ is monotonically increasing as a function of $\lambda\in(1,\infty)$ when $\psi\in(0,\pi/2)$ and monotonically decreasing when $\psi\in(-\pi/2,0)$, with the limits $\lim_{\lambda\searrow 1}\phi(\lambda,\lambda,\psi)=0$ and $\lim_{\lambda\nearrow\infty}\phi(\lambda,\lambda,\psi)=\psi$ in both cases.
\item[iii)]
   With $\lambda_1,\psi$ specified, 
   $\phi(\lambda_1,\lambda_2,\psi)$ is monotonically increasing as a function of $\lambda_2\in(0,\infty)$ when $\psi\in(0,\pi/2)$ and monotonically decreasing when $\psi\in(-\pi/2,0)$,
   with the limits in both cases given by
    \begin{align*}
        \lim_{\lambda_2\nearrow\infty}\phi(\lambda_1,\lambda_2,\psi) = \tan^{-1}\left(\frac{(\lambda_1-1)\tan \psi}{\lambda_1 + \tan^2 \psi}\right),\\
        \lim_{\lambda_2\searrow 0}\phi(\lambda_1,\lambda_2,\psi) = \tan^{-1}\left(\frac{(1-\lambda_1)\tan \psi}{1 + \lambda_1\tan^2 \psi}\right).
    \end{align*}
   \item[iv)] With $\lambda_1>1$ specified, we have that
   \begin{align}\label{eq:phisup}
\phi_{\rm sup}(\lambda_1):=\sup_{\lambda_2>0,\psi\in(-\pi/2,\pi/2)}\phi(\lambda_1,\lambda_2,\psi)&=\tan^{-1}\left(\frac12\left(\lambda_1^{\frac12}-\lambda_1^{-\frac12}\right)\right),
   \end{align}
   approached in the limit $\lambda_2\nearrow \infty$ with $\psi=\tan^{-1}(\lambda_1^{\frac{1}{2}})$, as well as the limit $\lambda_2\searrow 0$ with $\psi=-\cot^{-1}(\lambda_1^{\frac{1}{2}})$, and, similarly, 
   \begin{align}\label{eq:phiinf}
\phi_{\rm inf}(\lambda_1):=\inf_{\lambda_2>0,\psi\in(-\pi/2,\pi/2)}\phi(\lambda_1,\lambda_2,\psi)&=\tan^{-1}\left(\frac12\left(\lambda_1^{-\frac12}-\lambda_1^{\frac12}\right)\right),
   \end{align}
   approached in the limit $\lambda_2\nearrow \infty$ with $\psi=-\tan^{-1}(\lambda_1^{\frac{1}{2}})$, as well as the limit $\lambda_2\searrow 0$ with $\psi=\cot^{-1}(\lambda_1^{\frac{1}{2}})$.
\end{itemize}
\end{lemma}



\begin{figure}[h]
  \centering
  \includegraphics[width=0.35\textwidth]{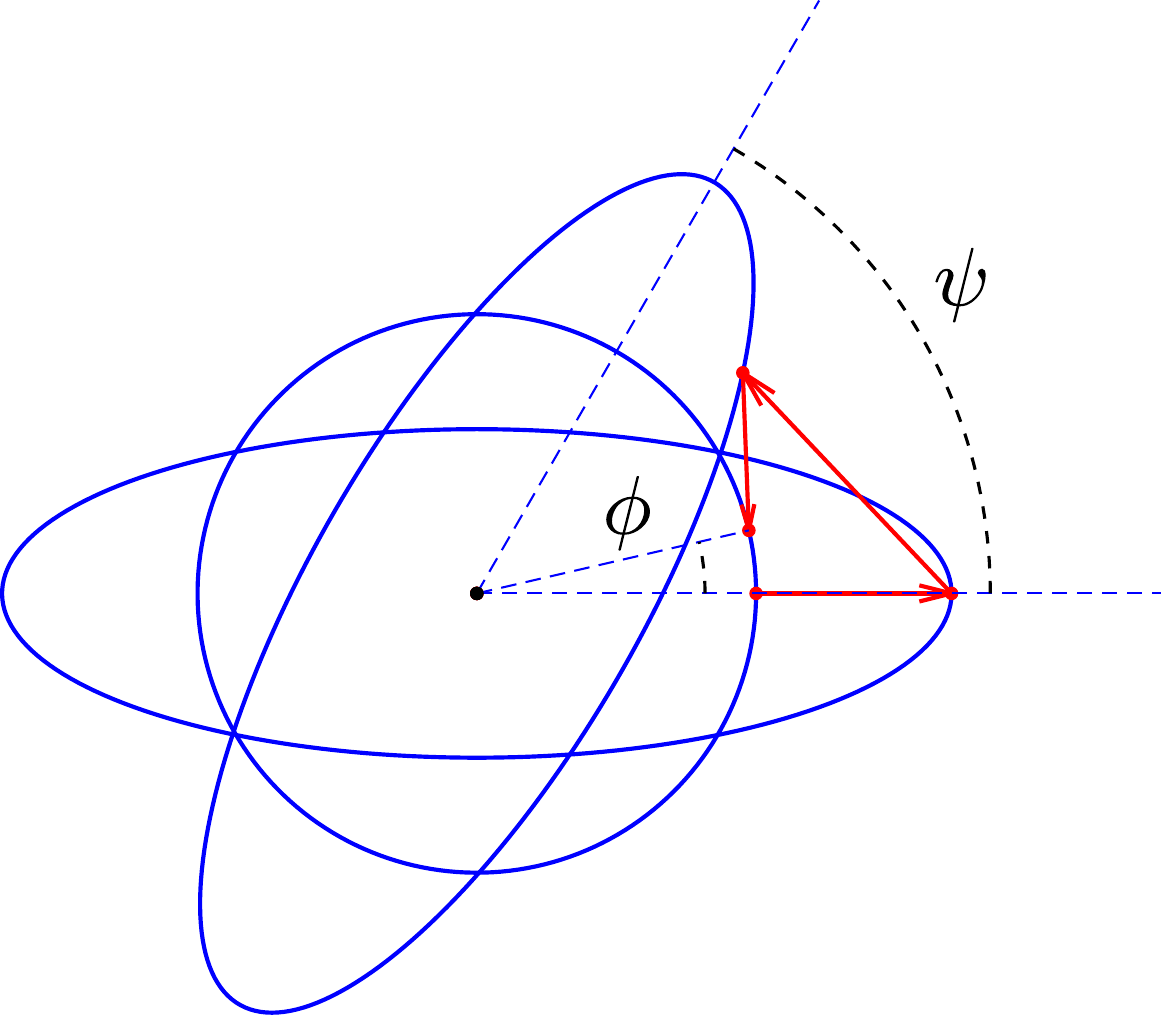}
    \caption{An illustration of the transitions $v\mapsto \Sigma_1^{1/2}v\mapsto \Phi_{21}\Sigma_1^{1/2}v\mapsto \Sigma_2^{-1/2}\Phi_{21}\Sigma_1^{1/2}v$.}
    \label{fig:fig2a}
    \end{figure}

\noindent{\em Proof}\\
i) It follows directly from \eqref{eq:R}, since $R_{12}=R_{21}^{-1}$.

ii) Figure \ref{fig:fig2a} exemplifies the successive transformations
for the case where $\lambda_1=\lambda_2=\lambda$.
In order to simplify typesetting we let $\Lambda:=\diag(\sqrt{\lambda},\sqrt{1/\lambda})$, and thus,
$\Sigma_1^{1/2}=\Lambda$ and $\Sigma_2=\Theta_\psi\Lambda\Theta^\top_\psi$. From \eqref{eq:key},
\begin{align*}
    \Theta_\phi&=\Theta_\psi \Lambda^{-1} \Theta_\psi^\top\Lambda^{-1}(\Lambda \Theta_\psi \Lambda^2 \Theta_\psi^\top \Lambda)^{1/2}\\
    \Rightarrow \Theta_{\phi-\psi}&=\Lambda^{-1} \Theta_\psi^\top\Lambda^{-1}(\Lambda \Theta_\psi \Lambda^2 \Theta_\psi^\top \Lambda)^{1/2} = M^{-1}(MM^\top)^{1/2},
\end{align*}
where 
\begin{equation}\label{eq:M}
M=\Lambda \Theta_\psi \Lambda =
\begin{bmatrix}
    \lambda\cos(\psi)&-\sin(\psi)\\\sin(\psi)&\frac{1}{\lambda}\cos(\psi)
\end{bmatrix}.
\end{equation}
Note that $\det(M)=1$, and so $\det(MM^\top)=1$ as well. We denote by $\mu$ and $1/\mu$ the eigenvalues of $MM^\top$, taking $\mu\geq 1$. 
Next, we determine $a,b$ so that
\[
(MM^\top)^{1/2}=aI+bMM^\top.
\]
To this end, we solve
$a+b\mu=\sqrt{\mu}$ and $a+b/\mu=\sqrt{1/\mu}$, giving that either $\mu=1$, or $a=b$. If $\lambda=1$, trivially $\Theta_\psi=I$. Else $\lambda>1$, and $\mu\neq 1$, and hence, 
\[
(MM^\top)^{1/2}=a(I+MM^\top),
\]
for $a=\sqrt{\mu}/(\mu+1)$.
 Therefore, $\Theta_{\phi-\psi}=a(M^{-1}+M^\top)$, and
 \begin{align}\label{eq:cosofdifference}
     \cos(\phi-\psi)&=\frac{1}{2}\trace\Theta_{\phi-\psi}=a\trace M=\frac{\sqrt{\mu}}{\mu+1}\cos(\psi)(\lambda+\frac{1}{\lambda}),\\
     \sin(\phi-\psi)&=-\frac{2\sqrt{\mu}}{\mu+1}\sin(\psi).\label{eq:sinofdifference}
 \end{align}
 Then
 $\tan(\phi-\psi) =-\frac{2\lambda}{\lambda^2+1}\tan\psi$,
and thus,
 \begin{align}\label{eq:tan2}
     \tan\phi &=\frac{\left(\lambda+\frac{1}{\lambda}-2\right)\tan\psi}{\lambda+\frac{1}{\lambda}+2\tan^2\psi}.
 \end{align}
Denote $\lambda+\frac{1}{\lambda}=:y\geq 2$, since we consider $\lambda\geq 1$, and note that $y$ monotonically increases as a function of $\lambda\in[1,\infty)$. The expression $\frac{y-2}{y+2\tan^2\psi}$ is monotonically increasing in $y$, and hence, in $\lambda$, readily giving that $\tan \phi(\lambda,\lambda,\psi)$ is monotonically increasing as a function of $\lambda\in(1,\infty)$ when $\psi\in(0,\pi/2)$ and monotonically decreasing when $\psi\in(-\pi/2,0)$, with $\tan \phi(1,1,\psi) = 0$ and $\lim_{\lambda\to\infty} \tan\phi(\lambda,\lambda,\psi)=\tan\psi$ in both cases.
We conclude that as $\lambda\to\infty$, $\phi\to\psi$ as claimed.

iii) Starting from a specified pair $(\lambda_1,\psi)$, with $\lambda_1>0$ and $\psi\in(-\pi/2,\pi/2)$, we proceed along identical lines to write, based on
\eqref{eq:key}, that once again
\[
\Theta_{\phi-\psi}=M^{-1}(MM^\top)^{1/2},
\]
though now,
\begin{equation}\label{eq:Mnew}
M=\Lambda_1 \Theta_\psi \Lambda_2 =
\begin{bmatrix}    \sqrt{\lambda_1\lambda_2}\cos(\psi)&-\sqrt{\frac{\lambda_1}{\lambda_2}}\sin(\psi)\\\sqrt{\frac{\lambda_2}{\lambda_1}}\sin(\psi)&\frac{1}{\sqrt{\lambda_1\lambda_2}}\cos(\psi)
\end{bmatrix},
\end{equation}
where
$\Lambda_i=\diag(\sqrt{\lambda_i},\sqrt{1/\lambda_i})$, with $i\in\{1,2\}$, and $\Sigma_1^{1/2}=\Lambda_1$ while $\Sigma_2^{1/2}=\Theta_\psi\Lambda_2\Theta_\psi^\top$. Looking at
\eqref{eq:cosofdifference},
mutatis mutandis,
\begin{align*}
    \cos(\phi-\psi)&=\frac{\sqrt{\mu}}{\mu+1}\cos(\psi)(\sqrt{\lambda_1\lambda_2}+\frac{1}{\sqrt{\lambda_1\lambda_2}})\\
    \sin(\phi-\psi)&=-\frac{\sqrt{\mu}}{\mu+1}\sin\psi \left(\sqrt{\frac{\lambda_1}{\lambda_2}}+\sqrt{\frac{\lambda_2}{\lambda_1}}\right),
\end{align*}
giving that $\tan(\phi-\psi)=-\frac{\lambda_1+\lambda_2}{\lambda_1\lambda_2+1}\tan\psi$, and, therefore, that
\begin{align*}
    \tan\phi(\lambda_1,\lambda_2,\psi)&=
    \frac{\left(\frac{\lambda_1\lambda_2+1}{\lambda_1+\lambda_2}-1\right)\tan\psi}{\frac{\lambda_1\lambda_2+1}{\lambda_1+\lambda_2}+\tan^2\psi}
    .
\end{align*}
For $\lambda_1>1$, the expression $y=\frac{\lambda_1\lambda_2+1}{\lambda_1+\lambda_2}$ is monotonically increasing in $\lambda_2$, and in turn, the expression $\frac{y-1}{y+\tan^2\psi}$ is monotonically increasing with $y$. Consequently, $\phi$ is monotonically increasing in $\lambda_2$ when $\tan \psi > 0$ and monotonically decreasing when $\tan \psi < 0$, with the limits in both cases given by
\begin{align}
    \lim_{\lambda_2\nearrow\infty}\tan\phi(\lambda_1,\lambda_2,\psi) = \frac{(\lambda_1-1)\tan \psi}{\lambda_1 + \tan^2 \psi}, \label{eq:phi_lim_infty}\\
    \lim_{\lambda_2\searrow 0}\tan\phi(\lambda_1,\lambda_2,\psi) = \frac{(1-\lambda_1)\tan \psi}{1 + \lambda_1\tan^2 \psi}. \label{eq:phi_lim_zero}
\end{align}
iv) The case $\psi=0$ is trivial and, thus, is omitted. From \eqref{eq:phi_lim_infty}, the supreme value of $\phi$ for a fixed $\lambda_1>1$ and any $\psi\in(0,\pi/2)$ is approached in the limit $\lambda_2\to\infty$ when $\psi=\tan^{-1}(\lambda_1^{\frac{1}{2}})$, and coincides with \eqref{eq:phisup}. Similarly, from \eqref{eq:phi_lim_zero}, the supreme value of $\phi$ for a given $\lambda_1>1$ and $\psi\in(-\pi/2,0)$ is also approached in the limit $\lambda_2\to 0$ when $\psi=-\cot^{-1}(\lambda_1^{\frac{1}{2}})$, and coincides with \eqref{eq:phisup}. On the other hand, the infimal value of $\phi$ for a given $\lambda_1>1$ and $\psi\in(0,\pi/2)$ is approached in the limit $\lambda_2\to\infty$ when $\psi=-\tan^{-1}(\lambda_1^{\frac{1}{2}})$, and coincides with \eqref{eq:phiinf}. Similarly, from \eqref{eq:phi_lim_zero}, the infimal value of $\phi$ for a given $\lambda_1>1$ and $\psi\in(-\pi/2,0)$ is also approached in the limit $\lambda_2\to 0$ when $\psi=\cot^{-1}(\lambda_1^{\frac{1}{2}})$, and coincides with \eqref{eq:phiinf}.
\hfill $\Box$


\section{Factorization into positive factors}\label{sec:4}

In light of the above analysis we outline a constructive approach to factor a given matrix with positive determinent into a product of positive factors. We begin with $2\times 2$ matrices.

\newcommand{\A}{A}
 Consider $\A \in\mathbb R^{2\times 2}$ where, as before, without loss of generality we assume that $\det \A=1$. 
Starting from the polar factorization, we write
 \begin{equation}\label{eq:polar}
 \A =\Theta_\chi (A^\top A)^{1/2}=\Theta_\chi \Theta D \Theta^\top =\Theta\Theta_\chi D\Theta^\top,
 \end{equation}
 so that $D=\diag(d,\,d^{-1})$ is diagonal with $d>1$, and $\Theta,\Theta_\chi$ are orthogonal matrices. We first observe that it suffices to factor
 $\Theta^\top \A  \Theta$ into a product
 \[
 \Phi_k\cdots \Phi_2\Phi_1
 \]
of positive factors, since then, $\A =\Theta\Phi_k\Theta^\top \cdots \Theta\Phi_2\Theta^\top\Theta\Phi_1\Theta$ readily gives a corresponding factorization for $\A $ with positive factors $\Theta\Phi_i\Theta^\top$ ($i\in\{1,\ldots,k\}$). Thus, without loss of generality, we henceforth assume that in \eqref{eq:polar},  $\Theta=I$.

We next observe that factoring $\A =\Theta_\chi D$ into a product of $k$ factors $\Phi_k\ldots\Phi_1$ is equivalent to factoring the corresponding orthogonal matrix $\Theta_\chi$ into $k+1$ positive factors with the first factor being $D^{-1}$, since then, $\Theta=\Phi_k\ldots\Phi_1D^{-1}$.
By appealing to Lemma \ref{lemma:3.1} we now have the following.

\begin{thm}
Let
$\A = \Theta_\chi D$
as above and
$\phi_{\rm sup}(d)=\tan^{-1}(d^{\frac{1}{2}}-d^{-\frac{1}{2}})$,
as in \eqref{eq:phimax}. The following statments hold:
\begin{itemize}
\item[a)] If $d< 1$ and $|\chi|<\phi_{\rm sup}(d)$, there exist a positive factorization $A=\Phi_2 \Phi_1$, with two symmetric factors.
\item[b)] If $d< 1$ and $|\chi|<\frac{\pi}{2}+\phi_{\rm sup}(d)$, there exist a positive factorization $A=\Phi_3\Phi_2 \Phi_1$, with three symmetric factors.

\item[c)] If $d< 1$ and $|\chi|\leq\pi$, there exist a positive factorization $A=\Phi_4\Phi_3\Phi_2 \Phi_1$, with four symmetric factors.
\item[d)] If  $d=1$ and $\chi=\pi$, i.e., $A=-I$, there exist a positive factorization $A=\Phi_5\Phi_4\Phi_3\Phi_2 \Phi_1$, with five symmetric factors.
\end{itemize}
Moreover, the conditions given in each case are also necessary.
\end{thm}

\noindent{\em Proof}
    a) From the discussion leading up to the theorem, we need to construct a factorization $\Theta_\chi=\Phi_3\Phi_2 \Phi_1$, starting from $\Phi_1=D^{-1}$ as the first factor.
    From Lemma \ref{lemma:3.1}.iv) the supremal (respectively, infimal) value for the accrued holonomy that can be obtained in the first step, i.e., by the composition of $\Phi_1=\Sigma_2^\frac12$ and $D=\Sigma_1^\frac12$, is $\phi_{\rm sup}(d)$ (respectively, $\phi_{\rm inf}(d)$). The third factor must be $\Phi_3=\left(
\Sigma_2^{1/2}
\Sigma_1
\Sigma_2^{1/2}
\right)^{-1/2}$, as dictated by the requirement that the product of the three factors is orthogonal, see \eqref{eq:R}. Thus, $\Theta_\chi$ can be constructed by a suitable choice of values for $\psi$ and $\lambda_2$, to specify $\Sigma_2$ as in \eqref{eq:Sigma12}. These are the only two parameters to be selected, and need to satisfy $\phi(d,\lambda_2,\psi)=\chi$, which is possible as readily follows from Lemma \ref{lemma:3.1} claims iii) and iv).

b,c) These cases are similar to a). We now need to factor $\Theta_\chi=\Phi_4\Phi_3\Phi_2 \Phi_1$, and $\Theta_\chi=\Phi_5\Phi_4\Phi_3\Phi_2 \Phi_1$, respectively, starting from $\Phi_1=D^{-1}$. The maximal holonomy accrued by choice of $\Phi_2$ in precisely the same manner is $\phi_{\rm sup}(d)$ or $\phi_{\rm inf}(d)$, depending on the required sign. In subsequent steps, i.e., selecting $\Phi_3$ and $\Phi_3$, $\Phi_4$, respectively, the maximal holonomy is $\pi/2$. The last factor in each product is dictated by the requirement that the product is an orthogonal matrix. Thus, as long as $|\chi|$ is within the specified bounds, the factors can be constructed mutatis mutandis.

d) In this case $D=I$, and the choice begins with the first factor $\Phi_1=\Sigma_1^{\frac12}$, diagonal, for a respective value $\lambda_1$. Then, a successive selection of $\psi_k,\lambda_k$, to give
\[
\Sigma_k= \Theta_{\psi_k} \begin{bmatrix}\lambda_k&0\\0&1/\lambda_k
\end{bmatrix}\Theta_{\psi_k}^\top,
\]
for $k\in\{2,3,4\}$ can provide holonomy $<\pi/2$ in each step. 
Clearly, three steps are needed to ensure that the total holonomy equals $\pi$ --two steps do not suffice. 
That a selection is feasible follows readily from statement ii) Lemma \eqref{lemma:3.1}. 
Specifically, select $\psi_\infty>\pi/3$ and $\lambda$ so that $\phi(\lambda,\lambda,\psi_\infty)=\pi/3$ and let $\Sigma_1=\diag(\lambda,1/\lambda)$, $\Sigma_{1+k}=\Theta_\psi^k \Sigma_1 \Theta_\psi^k$ for $k\in\{1,2,3\}$. These specify the first four factors that generate a holonomy of precisely $\pi$ [rad]. The last factor is dictated by the requirement that the product is an orthogonal matrix. That the holonomy is $\pi$ gives that this is precisely the negative of the identity.
\hfill $\Box$

Case a) in the theorem corresponds to the situation when $A$ has positive eigenvalues which can be readily verified as being equivalent to the condition provided.
Also, the factorizations in the theorem are evidently not unique. The freedom in selecting parameters can be used to ensure that the factors have the same conditioning number. 
    For instance, a five-factor factorization of $\A=-I$, with factors having (approximately) the same conditioning number can be obtained by choosing $\lambda=30$ and $\psi=70.3^o$. This choice gives that
  \begin{align*}
 \begin{bmatrix}
      -1 &0\\0&-1
  \end{bmatrix} = &\begin{bmatrix}
        1.58 &  -2.34\\
   -2.34 &   4.08
  \end{bmatrix}
  \begin{bmatrix}
      3.32&2.71\\2.71 & 2.52
  \end{bmatrix}\times\\
   &\times \begin{bmatrix}
      4.33&-2.35\\-2.35&1.50
  \end{bmatrix}
   \begin{bmatrix}
      .34 & .92\\.92 &5.50
  \end{bmatrix}
   \begin{bmatrix}
      5.48 &0\\0&0.18
  \end{bmatrix},
  \end{align*}
  with the factors rounded at two decimals; the eigenvalues of each factor are the same, and the condition number is $\lambda=30$.

\begin{remark}
    It is interesting to note that the perceived effect of rotation via a succession of gradient flows can be exemplified with working with playdough, and exerting compressive forces sequentially along different directions. It is a fact that one may have experienced while playing as a  child, that a mark on the playdough may rotate in position, without ever having rotated the bulk of the playdough in our hands. This has been the theme of the present work -- how the successive application of irrotational vector fields may conspire to produce a rotation.
\end{remark}

\section{Factorization of rotation matrices in $\mathbb R^n$}\label{sec:5}

We now consider the general problem on how to factor matrices of arbitrary size and positive determinant into a finite product of positive definite factors. We show that the problem reduces to factoring $2\times 2$ matrices in the same way.
We focus on orthogonal matrices, since polar factorization may provide one of the factors.
Thus, we consider an arbitrary $n\times n$ orthogonal matrix $V$ and seek a factorization of $V$ into a product of symmetric matrices. The factorization can be accomplished as explained below.

Since $V$ is orthogonal, and hence
\[
VV^T=V^TV=I,
\]
it is also normal; it possesses a complete set of orthogonal eigenvectors.
Thus, it admits a block-diagonal decomposition
\[
V= U\begin{bmatrix}
    \begin{matrix}
        \hspace{-3pt}\phantom{-}\cos(\theta_1) & \sin(\theta_1)\\
        \hspace{-3pt}-\sin(\theta_1) & \cos(\theta_1)
    \end{matrix} & 0 & \\
    0 & 
    \begin{matrix}
    \hspace{-7pt}\phantom{-}\cos(\theta_2) & \sin(\theta_2)\\
       \hspace{-7pt} -\sin(\theta_2) & \cos(\theta_2)&\\        0 & 0 & \ddots 
    \end{matrix}
\end{bmatrix}U^T
\]
into $2\times 2$ blocks when the size is even,
and
with one additional $1\times 1$ block with entry equal to $1$ when the size is odd. Note that if the multiplicity of the eigenvalues at $\pm 1$ is higher, then these can be accounted for with $2\times 2$ (diagonal) block entries.
All  $2\times 2$ block entries are rotation matrices, corresponding to rotations on the plane by angles $\theta_i$, ($i\in\{1,\ldots\})$. Possible reflections come in pairs that can be grouped into $2\times 2$ blocks that are diagonal with entries equal to $-1$. 
Also, $U$ is the orthogonal matrix that is made up of the eigenvectors of $V$.

Using the scheme explained in the previous section, each block can be factored into
a product of at most five $2\times 2$ symmetric matrices. In case the size of the matrix is odd, and there is an additional $1\times 1$, identity block, this factors similarly and trivially into five such identity blocks. Putting all these together, we have a factorization
\begin{align*}
V&=U M_5\cdots M_1 U^T\\
&=N_5\cdots N_1
\end{align*}
with symmetric factors $N_i:=UM_iU^T$, for $i\in\{1,\ldots,5\}$.

Thus, the steps we have described, based on the factorization in Section~\ref{sec:4}, readily provides a factorization of any matrix $\Phi$ with positive determinent into a product of no more than six factors. Using the scheme we provided, factorizations into more than six factors can be readily obtained, with factors having a better numerical condition (maximum ratio of corresponding eigenvalues).

We note that Ballantine \cite{ballantine1968aproducts,ballantine1970products} established that no more than five factors are needed in general. However, his argument is non-constructive. In this extreme case, when five factors are needed, a constructive method to obtain the factorization is still absent. The constructive approach outlined above for $n\times n$ matrices can be used to obtain six factors.

\section{Concluding remarks}\label{sec:remarks}

The purpose of this note is to explain and outline a systematic approach to factoring arbitrary matrices with positive determinant into products of symmetric factors. Such factorizations are useful in designing control protocols for collections of dynamical systems that follow identical dynamics as first sought in an influential paper by Roger Brockett \cite{brockett2007optimal} and completed by the present authors in \cite{abdelgalil2024collective}.  

The existence of such factorizations, with a guaranteed finite upper bound on the number of factors, were first established by Ballantine. In fact, Ballantine \cite{ballantine1968aproducts,ballantine1970products} provides conditions when a matrix can be given as a product of three and four factors, and provides a non-constructive proof that in general no more than five factors are needed for a matrix with positive determinant.

Thus, the present work attempts to provide a constructive approach to obtaining Ballantine factorizations.
Our approach is substantially different from that of Ballantine. Our approach is based on an explicit construction of the sought symmetric factors as maps that optimally relate suitably rotated Gaussian distributions, and was inspired by the Monge-Kantorovich transportation theory specialized to Gaussians. A theory that would allow factoring more general (nonlinear) diffeomorphisms as a composition of gradient flows is not available and it is a subject currently pursued by the authors.

\bibliographystyle{plain}
\bibliography{References}
\end{document}